\documentclass[11pt, final]{amsart}
\usepackage{amssymb}
\usepackage[mathscr]{eucal}

\theoremstyle{plain}
\newtheorem*{theorem}{Theorem}
\newtheorem*{corollary}{Corollary}

\theoremstyle{remark}
\newtheorem*{remark}{Remark}

\newcommand\C{\mathbb C}
\newcommand\R{\mathbb R}

\begin{document}
\title{${}^*$-semigroup endomorphisms of $B(H)$}
\author{LAJOS MOLN\'AR}
\address{Institute of Mathematics\\
         Lajos Kossuth University\\
         4010 Debrecen, P.O.Box 12, Hungary}
\email{molnarl@math.klte.hu}
\dedicatory{Dedicated to the memory of Professor B\'ela
Sz\H{o}kefalvi-Nagy}
\thanks{  This research was supported from the following sources:\\
          1) Hungarian National Foundation for Scientific Research
          (OTKA), Grant No. T--030082 F--019322,\\
          2) A grant from the Ministry of Education, Hungary, Reg.
          No. FKFP 0304/1997}
\date{\today}
\begin{abstract}
Let $H$ be a complex separable infinite dimensional Hilbert space.
We describe the form of
all *-semigroup endomorphisms $\phi$ of $B(H)$ which are uniformly
continuous on every commutative $C^*$-subalgebra.
In particular, we obtain that if $\phi$ satisfies $\phi(0)=0$, then
$\phi$ is additive.
\end{abstract}
\maketitle

\vskip 1cm
\section{Introduction}
In the paper \cite{Hak} the author studied the question
of additivity of *-semigroup isomorphisms between operator algebras.
The problem that when multiplicativity implies additivity was previously
investigated also from the purely algebraic
point of view \cite{Mar}. Using that result, the semigroup isomorphisms
between standard operator algebras were completely described in
\cite{Semrl}. The common characteristic of the mentioned investigations
is that they considered multiplicative bijections, or, as in \cite{Mar},
at least the surjectivity was assumed.
In our present paper we get rid of this very restrictive algebraic
condition and, instead, we impose a certain continuity assumption.
In contrast to the proofs of the previously mentioned results
which were based mainly on algebraic manipulations, our argument needs
much more functional analysis and this is the reason why we obtain our
result "only" for the full operator algebra over a separable
infinite dimensional Hilbert space
(as for the finite dimensional case, we should refer to the paper
\cite{JL} where all so-called non-degenerate multiplicative maps of a
matrix algebra were determined).

\section{Results}

The main result of the paper reads as follows.

\begin{theorem}
Let $H$ be a complex separable infinite dimensional Hilbert space.
Let $\phi: B(H) \to B(H)$ be a *-semigroup
endomorphism which is uniformly continuous on the commutative
$C^*$-subalgebras of $B(H)$.
Then $\phi$ can be written in the form
\[
\phi(A)=
\left[
\begin{matrix}
0       &   0      & \dots  &\dots&\dots&\dots&\dots&\dots&\dots \\
0       &   I      & 0&\dots&\dots&\dots&\dots&\dots&\dots\\
0       &   0      & A         & 0&\dots&\dots&\dots&\dots&\dots \\
0       &   0      & 0         & A&0    &\dots&\dots&\dots&\dots& \\
\vdots&\vdots&\vdots&\vdots&\ddots&\vdots&\vdots&\vdots&\vdots\\
0     & \dots&\dots &\dots &0     &{A^*}^{tr} &0&\dots&\dots \\
0     &\dots &\dots &\dots &\dots &0          &{A^*}^{tr}&0&\dots \\
\vdots&\vdots&\vdots&\vdots&\vdots&\vdots&\vdots&\ddots&\vdots\\
\end{matrix}
\right]
\quad (A\in B(H)).
\]
\end{theorem}

\begin{proof}
Clearly, $\phi$ sends projections to projections. So, $\phi(I)$ and
$\phi(0)$ are projections. Since
\[
\phi(I)\phi(A)=\phi(A)\phi(I)=\phi(A)
\]
and
\[
\phi(0)\phi(A)=\phi(A)\phi(0)=\phi(0),
\]
it is easy to see that
$\phi$ can be written in the form
\[
\phi(A)=
\left[
\begin{matrix}
0       &   0      & 0      \\
0       &   I      & 0      \\
0       &   0      & \phi'(A) \\
\end{matrix}
\right]
\qquad (A\in B(H))
\]
where $\phi'$ is a *-semigroup endomorphism of
$B(H)$ having the same continuity property as $\phi$
which sends 0 to 0 and maps $I$ into $I$.
Therefore, we can assume that our original map $\phi$
satisfies $\phi(0)=0$ and $\phi(I)=I$.
The main step of the proof which follows is to prove that $\phi$ is
orthoadditive on the set of all projections.
This means that $\phi(P+Q)=\phi(P)+\phi(Q)$ for any mutually
orthogonal projections $P,Q\in B(H)$. To see this, let $P\in B(H)$ be
an arbitrary projection and set $Q=I-P$. Consider the map
\[
\lambda \longmapsto \phi(e^{\lambda P})=\phi(Q+e^\lambda P)
\]
from $\R$ into the group of all invertible operators in $B(H)$. This is
a continuous one-parameter group and hence there is an operator  $T\in
B(H)$ such that
\[
\phi(Q+e^\lambda P)= e^{\lambda T} \qquad (\lambda \in \R)
\]
(see, for example, \cite[6.4.6 Proposition]{Pal}).
Since $\phi$ is *-preserving, we obtain that $e^{\lambda T}$ is
self-adjoint for all $\lambda\in \R$. This yields that $T$ is
also self-adjoint.
The norm and the spectral radius of any self-adjoint operator coincide.
So, from the continuity property of $\phi$ we deduce that
for every $\epsilon >0$ there is a $\delta >0$ such that
$\sup_{t\in \sigma(T)} |e^{\lambda t}-e^{\mu t}|<\epsilon$
if $|e^\lambda -e^\mu|<\delta$.
Therefore, the function $x \mapsto x^t$ is uniformly
continuous on the positive half-line for every $t\in \sigma(T)$. This
gives us that $\sigma(T)\subset \{ 0,1\}$. Consequently, $T$ is a
projection. Then we have $e^{\lambda T}=(I-T)+e^\lambda T$ and thus
\[
\phi(Q+e^\lambda P)=(I-T) +e^\lambda T \qquad (\lambda \in \R)
\]
or, equivalently,
\begin{equation}\label{E:speci}
\phi(Q+\epsilon P)=(I-T) +\epsilon T
\end{equation}
for every positive $\epsilon$. By the continuity property of $\phi$ we
obtain $\phi(Q)=I-T$.
As a particular case, we have
\begin{equation}\label{E:spricc}
\phi(\epsilon I)=\epsilon I \qquad (\epsilon >0).
\end{equation}
Indeed,
this follows from $\phi(0)=0$ and \eqref{E:speci}.
Therefore, $\phi$ is positive homogeneous,
and referring to \eqref{E:speci} again, if we divide by $\epsilon$ and
use the continuity property of $\phi$, then we arrive at $\phi(P)=T$.
Thus, we obtain
\begin{equation}\label{E:hamm}
\phi(P)+\phi(I-P)=I
\end{equation}
for every projection $P$ in $B(H)$. If
$P,Q$ are arbitrary projections with $PQ=QP=0$, then we infer from
the multiplicativity of $\phi$  and \eqref{E:hamm} that
\[
\phi(P)+\phi(Q)=\phi(P)\phi(P+Q)+\phi(I-P)\phi(P+Q)=
\]
\[
(\phi(P)+\phi(I-P))\phi(P+Q)=\phi(P+Q).
\]
Consequently, $\phi$ is orthoadditive on the set of all projections.

Since $\phi(I)=I$, it follows that $\phi$ sends unitaries to unitaries.
Consider the map $t \mapsto \phi(e^{it}I)$. Clearly, this is a
continuous one-parameter
unitary group. By Stone's theorem there is a self-adjoint operator $S\in
B(H)$ such that
\[
\phi(e^{it}I)=e^{itS} \qquad (t\in \R).
\]
Since $\phi(I)=I$, we have $e^{2\pi i S}=I$.
By spectral mapping theorem
this yields that the spectrum of $S$ consists of integers.
So, $S$ can be written in the form
$S=\sum_{k=-n}^{n} k P_k$,
where the $P_k$'s are pairwise orthogonal projections with
$\sum_{k=-n}^n P_k=I$ and $n$ is a suitable positive integer. We
compute
\[
\phi(e^{it}I)=
e^{ it \sum_{k=-n}^{n} k P_k}=
\]
\[
\prod_{k=-n}^{n} e^{ it k P_k}=
\prod_{k=-n}^{n} (I+ (e^{ it k}-1) P_k)=
\]
\[
I+\sum_{k=-n}^{n} (e^{ it k}-1) P_k=
\sum_{k=-n}^{n} e^{ it k} P_k
\qquad (t \in \R).
\]
Consequently, we have
\begin{equation}\label{E:kaja}
\phi(\lambda I)=
\sum_{k=-n}^{k=n} \lambda^k P_k \qquad (\lambda \in \C,
|\lambda|=1).
\end{equation}
From \eqref{E:spricc} and \eqref{E:kaja} we infer that
\begin{equation}\label{E:kaja1}
\phi(\lambda I)=
\sum_{k=-n}^{n} |\lambda| \biggl(\frac{\lambda}{|\lambda|}\biggr)^k
P_k
\qquad (\lambda \in \C, \lambda\neq 0).
\end{equation}
We know that
$\phi(\lambda I)$ commutes with $\phi(A)$ for every $\lambda\in \C$.
Thus, for any $A\in B(H)$ we have
\[
\sum_{k,l=-n}^n \lambda^k P_k\phi(A)P_l=
\sum_{k,l=-n}^n \lambda^k P_k P_k\phi(A)P_l=
\]
\[
\phi(\lambda I)\phi(A)=
\phi(A)\phi(\lambda I)=
\]
\[
\sum_{k,l=-n}^n P_k \phi(A)P_l P_l \lambda^l P_l=
\sum_{k,l=-n}^n \lambda^l P_k \phi(A)P_l
\]
for every $\lambda \in \C$ of modulus 1. This implies that
$P_k\phi(A)P_l=0$ if $k\neq l$. Consequently,
$\phi$ can be written in the form
\[
\phi(A)=
\sum_{k=-n}^{k=n} P_k \phi(A) P_k \qquad (A\in B(H))
\]
or, in another way,
\[
\phi(A)=
\left[
\begin{matrix}
\phi_{-n}(A) &   0           & \dots & 0 \\
0            &\phi_{-n+1}(A) & \dots & 0 \\
\vdots       & \vdots &  \ddots & \vdots \\
0            &   \dots  & \dots & \phi_n(A) \\
\end{matrix}
\right]
\qquad (A\in B(H)).
\]
Here, every $\phi_k$ $(k=-n, \ldots, n)$ is
a *-semigroup endomorphism of $B(H)$ which is uniformly continuous on
the commutative $C^*$-subalgebras.

Every orthoadditive projection valued measure on the set of
all projections in $B(H)$ can be extended to a linear map on $B(H)$.
This is a particular case of the solution of the
Mackey-Gleason problem obtained in \cite{BW}.
Since this extension is linear and sends projections to projections,
it is a standard argument to verify that this is in fact a
Jordan *-homomorphism (see, for example, the proof of
\cite[Theorem 2]{MolStud}).
A linear map $J$ between *-algebras $\mathcal A$ and $\mathcal B$
is called a Jordan *-homomorphism if it satisfies
\[
J(x)^2=J(x^2), \quad J(x)^*=J(x^*) \quad (x\in \mathcal A).
\]
So, we have Jordan *-homomorphisms $\psi_{-n}, \ldots, \psi_n$ of
$B(H)$ such that
\[
\phi(P)=
\left[
\begin{matrix}
\psi_{-n}(P) &   0           & \dots & 0 \\
0            &\psi_{-n+1}(P) & \dots & 0 \\
\vdots       & \vdots &  \ddots & \vdots \\
0            &   \dots  & \dots & \psi_n(P) \\
\end{matrix}
\right]
\]
for every projection $P\in B(H)$.

Let $R_1, \ldots, R_m$ be pairwise orthogonal projections whose sum is
$I$ and pick nonzero scalars $\lambda_1, \ldots, \lambda_m\in \C$. Using
the orthoadditivity of $\phi$, for any $k=-n, \ldots, n$ we compute
\[
\phi_k(\lambda_1R_1 +\ldots +\lambda_mR_m)=
\phi_k(\lambda_1R_1 +\ldots +\lambda_mR_m)
\phi_k(R_1 +\ldots +R_m)=
\]
\[
\phi_k(\lambda_1R_1 +\ldots +\lambda_mR_m)
(\phi_k(R_1) +\ldots +\phi_k(R_m))=
\]
\[
\phi_k(\lambda_1R_1) +\ldots +\phi_k(\lambda_mR_m)=
\]
\[
\phi_k(\lambda_1I)\phi_k(R_1) +\ldots +\phi_k(\lambda_mI)\phi_k(R_m)=
\]
\[
|\lambda_1| \biggl(\frac{\lambda_1}{|\lambda_1|}\biggr)^k\phi_k(R_1) +
\ldots +
|\lambda_m| \biggl(\frac{\lambda_m}{|\lambda_m|}\biggr)^k\phi_k(R_m)=
\]
\[
|\lambda_1| \biggl(\frac{\lambda_1}{|\lambda_1|}\biggr)^k\psi_k(R_1) +
\ldots +
|\lambda_m| \biggl(\frac{\lambda_m}{|\lambda_m|}\biggr)^k\psi_k(R_m)=
\]
\[
\psi_k\biggl(|\lambda_1|
\biggl(\frac{\lambda_1}{|\lambda_1|}\biggr)^kR_1 + \ldots +
|\lambda_m| \biggl(\frac{\lambda_m}{|\lambda_m|}\biggr)^kR_m\biggr).
\]
Using the continuity property of $\phi_k$, the automatic continuity of
Jordan
*-homomorphisms between $C^*$-algebras and the spectral theorem of
normal operators, we deduce that
\begin{equation}\label{E:multi}
\phi_k(N)=\psi_k(|N|(N|N|^{-1})^k)
\end{equation}
holds for every invertible normal operator $N\in B(H)$
(note that $N$ and the range of its spectral measure generate a
commutative $C^*$-subalgebra).

Every Jordan *-homomorphism of
$B(H)$ is the direct sum of a *-homomorphism and a *-antihomomorphism
(see \cite[Theorem 3.3]{Stor}). Let $\psi_k^h$ denote the *-homomorphic
and let $\psi_k^a$ denote the *-antihomomorphic part of $\psi_k$.
Let $N,M\in B(H)$ be invertible normal operators whose product is
also normal.
By the multiplicativity of $\phi_k$ and \eqref{E:multi} we have
\[
\psi_k(|NM|(NM|NM|^{-1})^k)=
\phi_k(NM)=\phi_k(N)\phi_k(M)=
\]
\[
\psi_k(|N|(N|N|^{-1})^k)
\psi_k(|M|(M|M|^{-1})^k).
\]
This implies that
\[
\psi_k^h(|NM|(NM|NM|^{-1})^k)=
\]
\begin{equation}\label{E:normszor1}
\psi_k^h(|N|(N|N|^{-1})^k)
\psi_k^h(|M|(M|M|^{-1})^k)=
\end{equation}
\[
\psi_k^h(|N|(N|N|^{-1})^k |M|(M|M|^{-1})^k)
\]
and
\[
\psi_k^a(|NM|(NM|NM|^{-1})^k)=
\]
\begin{equation}\label{E:normszor2}
\psi_k^a(|N|(N|N|^{-1})^k)
\psi_k^a(|M|(M|M|^{-1})^k)=
\end{equation}
\[
\psi_k^a(|M|(M|M|^{-1})^k|N|(N|N|^{-1})^k).
\]
Any *-homomorphism or *-antihomomorphism of $B(H)$ is either
injective or identically 0 which follows from the form of
representations of $B(H)$ on separable Hilbert spaces (see
\cite[10.4.14. Corollary]{KR}). Now, taking \eqref{E:normszor1} and
\eqref{E:normszor2} into account, one can verify that
the only values of $k$ for which $\psi_k$ can be nonzero are $-1$ and
$1$.
Moreover, because of the same reasons, for $k=1$ we have $\psi_1^a=0$
and for $k=-1$ we have $\psi_{-1}^h=0$. Observe that
$|N|(N|N|^{-1})^{-1}=N^*$.
Therefore, $\phi$ can be written in the form
\begin{equation}\label{E:felcsopvan}
\phi(N)=
\left[
\begin{matrix}
\psi'(N) & 0\\
0      & \psi''(N^*) \\
\end{matrix}
\right]
\end{equation}
for every invertible normal operator $N$, where $\psi'$ is a
*-endomorphism and $\psi''$ is a *-antiendomorphism of $B(H)$.
By continuity and spectral theorem we clearly have \eqref{E:felcsopvan}
for every normal operator in $B(H)$.
Define
\[
\psi(A)=
\left[
\begin{matrix}
\psi'(A) & 0\\
0      & \psi''(A^*) \\
\end{matrix}
\right]
\]
for every $A\in B(H)$. Clearly, $\psi$ is an additive
*-semigroup endomorphism
of $B(H)$ (it is not linear unless $\psi''$ is missing).

It is easy to see that every rank-one operator is the product
of at most three (rank-one) normal operators. This gives us
that $\phi(A)=\psi(A)$ for every rank-one operator $A\in B(H)$.
Now, let $A\in B(H)$ be arbitrary.
Pick rank-one projections $P,Q\in B(H)$. Since $PAQ$ is of rank at
most 1, we compute
\begin{equation}\label{E:rekau}
\phi(P)\phi(A)\phi(Q)=
\phi(PAQ)=
\psi(PAQ)=
\end{equation}
\[
\psi(P)\psi(A)\psi(Q)=
\phi(P)\psi(A)\phi(Q).
\]
Since every *-endomorphism of $B(H)$ is normal, that is, weakly
continuous on the bounded subsets of $B(H)$ (see \cite[10.4.14.
Corollary]{KR}), it follows from \eqref{E:felcsopvan}
that for every maximal family
$(P_n)_n$ of pairwise orthogonal rank-one
projections we have $\sum_n \phi(P_n)=I$. Therefore, we infer from
\eqref{E:rekau} that $\phi(A)=\psi(A)$ for every $A\in B(H)$.
Finally, in order to get the explicite form of $\phi$
we refer once again to the form of linear *-endomorphisms of $B(H)$
appearing in
\cite[10.4.14. Corollary]{KR} and note that one could get
the form of linear *-antiendomorphisms of $B(H)$ in a similar way.
\end{proof}

As for the additivity of *-semigroup endomorphisms of $B(H)$ we obtain
the following result.

\begin{corollary}
Let $H$ be a complex separable infinite dimensional Hilbert space.
Let $\phi: B(H) \to B(H)$ be a *-semigroup
endomorphism which is uniformly continuous on the commutative
$C^*$-subalgebras of $B(H)$. If $\phi(0)=0$, then $\phi$ is
additive.
\end{corollary}

\begin{remark}
In the paper \cite{HakPAMS} the author studied the problem of additivity
of Jordan *-isomorphisms between operator algebras under the same
continuity
assumption that we have used in our paper. It was conjectured by S.
Sakai and proved in \cite{HakJapan} that this condition is
in fact redundant. Therefore, it would be interesting to study our
problem without assuming any kind of continuity. However, taking
into account the results in \cite{JL},
one can get evidence to see that our problem in that case is much more
complicated.

Finally, we note that it would also be interesting to consider more
general algebras like von Neumann algebras in place of $B(H)$ and thus
generalize our result for that case.
\end{remark}

\end{document}